\documentclass[12pt]{article}
\usepackage{amsfonts,amssymb,amsthm}
\newcommand{\vs}{\vspace{5pt}}
\newcommand{\n}{\noindent}
\newcommand{\q}{\quad}
\newcommand{\qq}{\qquad}
\newcommand{\qand}{\q\mbox{ and }\q}
\newcommand{\R}{\mathbb{R}}
\newcommand{\ol}{\overline}
\newcommand{\al}{\alpha}
\newcommand{\be}{\beta}
\newcommand{\ga}{\gamma}
\newcommand{\Ga}{\Gamma}
\newcommand{\Th}{\Theta}
\newcommand{\om}{\omega}
\newcommand{\si}{\sigma}
\newcommand{\oet}{\ol\eta}
\newcommand{\oom}{\ol\om}
\newcommand{\osi}{\ol\si}
\newcommand{\op}{\oplus}
\newcommand{\ot}{\otimes}
\newcommand{\cC}{\mathcal{C}}
\newcommand{\cg}{\mathcal{G}}
\newcommand{\co}{\mathcal{O}}
\newcommand{\fc}{\mathfrak{c}}
\newcommand{\g}{\mathfrak{g}}
\newcommand{\fh}{\mathfrak{h}}
\newcommand{\ft}{\mathfrak{t}}
\renewcommand{\ge}{\geqslant}
\renewcommand{\le}{\leqslant}
\newcommand{\we}{\wedge}
\newcommand{\We}{\hbox{\small$\bigwedge$}}
\newcommand{\Abel}{\mathrm{Abel}}
\newcommand{\Kur}{\mathrm{Kur}}
\newcommand{\bC}{\mathbb{C}}

\newcommand{\bs}{\backslash}
\newcommand{\lb}{\left[}
\newcommand{\rb}{\right]}
\renewcommand{\ll}{\langle\!\langle}
\newcommand{\rr}{\rangle\!\rangle}
\renewcommand{\oe}{\ol e}
\newcommand{\half}{\hbox{$\frac12$}}
\newcommand{\bt}{\mathbf{t}}
\newcommand{\bmu}{\hbox{\boldmath$\mu$}}
\newcommand{\bnu}{\hbox{\boldmath$\nu$}}
\newcommand{\bph}{\hbox{\boldmath$\phi$}}
\newcommand{\bPh}{\hbox{\boldmath$\Phi$}}
\newcommand{\dbar}{\overline\partial}
\newcommand{\eps}{\epsilon}
\newcommand{\suml}[1]{\hbox{$\sum\limits_{#1}\,$}}
\newcommand{\oper}[2]{\newcommand{#1}{\mathop{\mathrm{#2}}\nolimits} }
\oper{\im}{Im}
\oper{\Aut}{Aut}
\oper{\image}{Im}

\newtheorem{corollary}{Corollary}
\newtheorem{proposition}{Proposition}
\newtheorem{theorem}{Theorem}
\newtheorem{definition}{Definition}

\newtheorem{lemma}{Lemma}
\newcommand{\bproof}{\n{\it Proof: }}
\newcommand{\eproof}{\q QED\vs}
\theoremstyle{remark}

\newtheorem{exa}{Example}

\begin{document}

\title{Deformation of 2-Step Nilmanifolds\\
with Abelian Complex Structures} 

\author{C.~McLaughlin$^{1,2}$ \and H.~Pedersen$^3$ \and Y.S.~Poon$^2$ \and
S.~Salamon$^3$}

\date{}
\maketitle

\footnotetext[1]{some results herein were reported in \cite{Colin} to fulfill
a degree requirement}
\footnotetext[2]{partially supported by NSF DMS-0204002}
\footnotetext[3]{partially supported by EC contract HPRN-CT-2000-00101}

\small\n{\bf Abstract.} We develop deformation theory for abelian invariant
complex structures on a nilmanifold, and prove that in this case the
invariance property is preserved by the Kuranishi process. A purely algebraic
condition characterizes the deformations leading again to abelian structures,
and we prove that such deformations are unobstructed. Various examples
illustrate the resulting theory, and the behavior possible in 3 complex
dimensions.\vs

\n{AMS Subject Classification: 32G05; 53C15, 53C56, 57S25, 17B30}
\normalsize

\section{Introduction}

In this paper, we study complex structures associated to compact quotients of
nilpotent groups. These manifolds are called nilmanifolds, and an
investigation of the special class of Kodaira manifolds was completed in
\cite{GPP}. The present paper opens this discussion to a wider class of
nilmanifolds.

A left-invariant complex structure on a Lie group is said to be \emph{abelian}
if the complex space of (1,0)-vectors is an abelian algebra with respect to
Lie bracket. This definition only makes sense in the algebraic setting, and in
this context it is of particular interest to know to what extent a study of
\emph{invariant} complex structures on nilmanifolds captures the general
situation.

There is a total of six 6-dimensional 2-step nilpotent groups admitting
abelian complex structures \cite{simon}. If $R^n$ denotes a $n$-dimensional
abelian group and $H_{2n+1}$ a $(2n+1)$-dimensional Heisenberg group, then the
2-step nilpotent groups with abelian complex structures are $R^6$, $H_5\times
R^1$, $H_3\times R^3$, $H_3\times H_3$, the Iwasawa group $W_6$ and one
additional group which we denote by $P_6$. These groups are the 6-dimensional
instances of the respective series: $R^{2n}$, $H_{2n+1}\times R^{2m-1}$,
$H_{2n+1}\times H_{2m+1}$, $W_{4N+2}$ and $P_{4N+2}$. For example, $W_{4N+2}$
is the real group underlying a generalized complex Heisenberg group.

The compact quotients of $R^{2n}$ are complex tori, and their deformation and
moduli are well studied \cite{B-L}. A detailed account of the moduli space of
complex structures of a special compact quotient of $H_{2n+1}\times R$ was
recently given in \cite{GPP}, and in \cite{KS} a somewhat different method is
directed towards the study of $W_6$. The present paper helps to unite these
two approaches, particularly via the examples in \S6.

This paper presents a general approach to computing the deformations of 2-step
nilmanifolds with abelian complex structures. To analyze data on deformation
theory, our first step is to identify the Dolbeault cohomology of a 2-step
nilmanifold with abelian complex structure with the appropriate Lie algebra
cohomology (Theorem~1). General results of this nature were proved in
\cite{CF}, though we shall need our own explicit description of this
identification. The second step is to extend this identification to the
determination of harmonic representatives for Dolbeault cohomology with
coefficients in the tangent sheaf (Theorem~\ref{harmonic rep}). 

We derive the first main result of this paper in \S4.3:

\vs\n{\bf Theorem}\it\q Let $G$ be a 2-step nilpotent Lie group with
co-compact subgroup $\Ga$. Then any abelian invariant complex structure on
$X=\Ga\bs G$ has a locally complete family of deformations consisting entirely
of invariant complex structures.\rm\vs

This is proved by showing that an application of Kuranishi's method does not
take one outside the subspace of invariant tensors. The theorem implies that
any deformation of an abelian invariant complex structure is necessarily
equivalent to an invariant one, at least of the deformation is sufficiently
small.

Given this result, it makes sense to ask under what conditions the deformed
invariant structures remain abelian. Indeed, our techniques enable us to prove
that deformations preserving the abelian property are always unobstructed and
faithfully represented at the infinitesimal level:

\vs\n{\bf Theorem}\it\q On a 2-step nilmanifold $X$ with an abelian complex
structure, a vector in the virtual parameter space $H^1(X,\Th_X)$ is
integrable to a 1-parameter family of abelian complex structures if and only
if it lies in a linear subspace defining the abelian condition
infinitesimally. \rm\vs

The `only if' part is obvious, but the force of this result is the backwards
implication. Once we convert to Lie algebra cohomology, it reduces the
constraints of abelian deformations to purely algebraic equations that we
introduce in \S5 and collectively call `Condition~A'. Using this, one may
carry out an effective computation in terms of structural constants of the
nilpotent groups in question.

Further analysis of abelian deformations yields a characterization of the
Kodaira manifolds (defined in \S2.1) as those corresponding to a Lie algebra
with 1-dimensional center and for which all infinitesimal parameters are
integrable and lead to abelian deformations. The precise statement is
Theorem~\ref{fully} in \S5.1.

In the final section, we compute the relevant cohomology dimensions for a
number of nil 6-manifolds, each equipped with a natural abelian complex
structure. In so doing, we are able to compare the techniques of this article
with those of \cite{simon}, but we emphasize a complication that arises from a
choice of complex structure with added symmetry. The first theorem above
allows one to dispense with the Kuranishi method in the explicit construction
of parameter spaces, and replace it with a more direct calculation involving
invariant differential forms. This we do in Example~8, after having first
illustrated the power of the second theorem above in estimating the dimension
of a potential moduli space.

\section{Abelian complex structures}

Suppose that a Lie algebra $\g$ admits an endomorphism $J$ such that
\begin{equation}\label{ab} J\circ J=-1\qand[JA,JB]=[A,B]\end{equation} for all
$A,B$ in $\g$. It can be extended by left-translation to an endomorphism of
the entire tangent bundle of $G$. Then $J$ defines an invariant almost complex
structure on the group $G$ which is integrable, since (\ref{ab}) implies the
vanishing of the Nijenhuis tensor. A complex structure satisfying (\ref{ab})
is called {\it abelian}, and the identity also implies that the center $\fc$
is $J$-invariant.

Now assume that the Lie algebra is 2-step nilpotent. In particular, the first
derived algebra is contained in the center. Taking the quotient of the algebra
$\g$ with respect to the center, we obtain an abelian algebra $\ft$. When the
complex structure is abelian, it induces a complex structure on $\ft$. The
identities \[\g^{(1,0)}=\ft^{(1,0)}\op\fc^{(1,0)}\qand \g^{(0,1)}=\ft^{(0,1)}
\op \fc^{(0,1)}\] concerning type (1,0) and (0,1) vectors are therefore valid
at the level of vector spaces.

Let $\{X_i,JX_i: 1\le i\le n\}$ be a real basis for $\ft$ and $\{Z_\al,
JZ_\al: n+1\le\al\le n+m \}$ a real basis for $\fc$. A basis of $(1,0)$
vectors for the complex tangent bundle of $G$ is composed of the elements
\begin{equation}\label{TW} T_j=\half(X_j-iJX_j)\qand W_\al=\half
(Z_\al-iJZ_\al).\end{equation} The complex structural constants $E_{kj}^\al$
and $F_{kj}^\al$ are defined by \begin{equation}\label{eq:structural} \lb\ol
T_k,T_j\rb=\suml\al E_{kj}^\al W_\al+\suml\al F_{kj}^\al\ol W_\al,
\end{equation} and satisfy \[ F_{kj}^\al=-\ol E_{jk}^\al.\]

We continue to use Roman indices in the range $1,\ldots,n$ and Greek indices
for $n+1,\ldots,n+m$. Let $\om^j$ be the left-invariant (1,0)-forms dual to
the vectors $T_j$ for $1\le j\le n$ and annihilating the $W_\al$. They span
the space $\ft^{*(1,0)}$. Similarly, there are left-invariant (1,0)-forms
$\om^\al$ dual to $W_\al$, annihilating the $T_j$. The dual form of the
structural equation (\ref{eq:structural}) is \begin{equation}
\label{differential} d\om^\al=\suml{i,j}E_{ji}^\al \om^i\we\oom^j.
\end{equation} The forms $\om^i$ are all exact, and $\partial\om^\al=0$. Thus,

\begin{lemma}\label{holomorphic forms}The forms $\oom^1,\ldots,\oom^n,
\oom^{n+1},\ldots,\oom^{n+m}$ are all $\dbar$-closed.\end{lemma}\medbreak

Now suppose that there exists a discrete subgroup $\Ga$ of $G$ such that the
left quotient space $\Ga\bs G$ is compact. The resulting quotient is called a
\emph{nilmanifold}. Since the complex structure $J$ is left-invariant, it
descends to a complex structure on $X=\Ga\bs G$.  Such a discrete subgroup
always exists if there is a basis such that the real structural constants are
rational \cite{Malcev}.

Later in this paper, we shall study the deformation theory on such compact
complex manifolds. At an appropriate juncture (in \S4.1), we shall find it
convenient to introduce an invariant Hermitian metric on $X$. We shall choose
such a metric so that $\{ X_j, JX_j, Z_\al, JZ_\al\}$ forms a Hermitian
frame. First we describe some simple examples of nilmanifolds and complex
structures.

\subsection{Kodaira manifolds and other examples}

Our first example of a compact nilmanifold with an abelian complex structure
is a Kodaira manifold, a generalization of a Kodaira surface. We proceed to
list algebraic constructions of this and similar examples.

\begin{exa} On the vector space $\R^{2n+2}$ with basis $\{X_j,Y_j,Z,A\}$,
define a Lie algebra by setting \[[X_j, Y_j]=-[Y_j,X_j]=Z,\qq 1\le j\le n,\]
and declaring all other brackets to be zero. This turns the vector space into
the direct sum $\g=\fh_{2n+1}\op\ft_1$ of the Heisenberg algebra and the
1-dimensional algebra.

We define an almost complex structure $J$ on the Lie algebra $\g$ by means of
the equations \begin{equation}\label{J}\begin{array}{c} JX_j=Y_j,\q
JZ=A,\end{array}\end{equation} so that the equations $JY_j=-X_j$ and $JA=-Z$
are also understood. The endomorphism $J$ defines an abelian complex structure
on $\g$, and therefore on the manifold $H_{2n+1}\times R$ and a compact
quotient thereof. In the niotation (\ref{TW}), its complex structure equation
is \begin{equation} \label{hei-structure} [\ol T_j,T_j]=-\half i(W+ \ol
W). \end{equation} The moduli problem of the compact quotient of such complex
manifolds was studied extensively in \cite{GPP}.\end{exa}

We next describe a more general extension of the Heisenberg group, and then a
product example.

\begin{exa} Let $H_{2n+1}\times R^{2m+1}$ be the product of a real Heisenberg
group and an abelian Lie group with dimensions as specified. Let $\{X_j, Y_j,
Z\}$ be a basis for $\fh_{2n+1}$, and let $\{Z_0, Z_{2k-1}, Z_{2k}\}$ (with
$1\le k\le m$) be a basis for $R^{2m+1}$. The non-zero structural constants
are determined by the single set of equations \[ [X_j, Y_j]=Z.\] An abelian
complex structure is defined by setting \begin{equation}\label{JJ} JX_j=Y_j,\q
JZ=Z_0,\q JZ_{2k-1}=Z_{2k},\end{equation} in analogy to (\ref{J}).\end{exa}

\begin{exa}\label{HH} The product $H_{2n+1}\!\times\!H_{2m+1}$ is a 2-step
nilpotent group with 2-dimensional center. Let $\{X_j, Y_j, Z_1, X_k, Y_k,
Z_2\}$ (with $1\le j \le n, 1\le k\le m$) be a basis for its Lie algebra
$\fh_{2n+1}\op \fh_{2m+1}$. The non-zero structural constants are \[[X_j,
Y_j]=Z_1\qand [X_k,Y_k]=Z_2.\] Define an almost complex structure $J$ on this
space by \begin{equation} \label{JJJ} JX_j=Y_j,\q JX_k=Y_k,\q JZ_1=Z_2;
\end{equation} once again this defines an abelian complex structure.\end{exa}

We describe the remaining classes of examples in 6 dimensions for simplicity.

\begin{exa}\label{W6} The group structure of $W_6$ underlies that of the
complex Heisenberg group. On the algebra level, the structural equations of
$W_6$ are \begin{equation}\label{W6se} \lb X_1,X_3\rb =-\half Z_1,\q\lb
X_1,X_4\rb=-\half Z_2,\q\lb X_2,X_3\rb=-\half Z_2,\q \lb X_2,X_4\rb=\half
Z_1.\end{equation} An abelian complex structure is defined by \begin{equation}
\label{J1} JX_1=X_2,\q JX_3=-X_4,\q JZ_1=-Z_2,\end{equation} and this is
denoted $J_1$ in \cite{KS}. Beware that $J$ is not the standard bi-invariant
complex structure $J_0$ that makes $W_6$ a complex Lie group; indeed
\[J_0\lb A,B\rb=\lb J_0A,B\rb,\qq A,B\in\g,\] and so $J_0$ is definitely not
abelian. Nonetheless, both $J_0$ and $J_1$ induce the same orientation on
$W^6$.\end{exa}

\begin{exa}\label{P6} The structural equations for $P_6$ are given by \[[X_1,
X_2]=-\half Z_1,\q [X_1,X_4]=-\half Z_2,\q [X_2,X_3]=-\half Z_2,\] and
correspond to a degeneration of (\ref{W6se}). An abelian complex structure $J$
on $P_6$ is defined by \begin{equation} \label{JP} JX_1=X_2,\q JX_3=-X_4, \q
JZ_1=-Z_2.\end{equation} The associated Hermitian manifold $\Ga\bs P_6$ was
studied in \cite[\S5]{AGS}, where it was shown that $J$ is only one of a
\emph{finite} number of complex structures compatible with a fixed Riemannian
metric.\end{exa}

\section{Cohomology theory}

In order to perform deformation theory on the compact complex nilmanifold $X$,
we need to calculate cohomology with coefficients in the holomorphic tangent
sheaf. We achieve this by identifying Dolbeault and Lie algebra cohomology, in
the spirit of \cite{Nom}.

\subsection{Lie algebra cohomology}

With respect to a complex structure $J$, the complexified Lie algebra has a
type decomposition. We may write \[\g_\bC=\g^{1,0}\op\g^{0,1},\q
\ft_\bC=\ft^{1,0}\op\ft^{0,1},\q \fc_\bC=\fc^{1,0}\op\fc^{0,1}.\] These are
all spaces of left-invariant vectors on $G$. The definitions are extended to
invariant $(p,q)$-forms in the standard way. For instance, $\We^k
\g_{\bC}^{*(0,1)}=\g^{*(0,k)}$ is the space of $G$-invariant $(0,k)$-forms.

Motivated by the property of the Chern connection on holomorphic tangent
bundles \cite{Gaud}, we define a linear operator $\dbar$ on $(0,1)$-vectors as
follows. For any $(1,0)$-vector $V$ and $(0,1)$-vector $\bar{U}$, set
\[\dbar_{\bar U}V:=[\bar U,V]^{1,0}.\] We obtain a linear map \[\dbar:
\g^{1,0}\to\g^{*(0,1)}\ot\g^{1,0}.\] In view of (\ref{eq:structural}), 
\begin{equation}\label{eq:dbar} \dbar_{\ol{T}_k}T_j=[\ol{T}_k,T_j]^{1,0}=
\suml\al E_{kj}^\al W_\al,\end{equation} whence \[\dbar T_j=\suml{k,\al}
E_{kj}^\al\ol{\om}^k\ot W_\al\qand\dbar W_\al=0.\]

Extend this definition to a linear map on $\g^{*(0,k)}\ot\g^{1,0}$ by setting
\[\dbar(\oom\ot V)=\dbar\oom \ot V+(-1)^k\oom\we\dbar V,\] where $V\in
\g^{*(0,k)}$ and $V\in\g^{1,0}$. For instance, any element $\bmu$ in
$\g^{*(0,1)}\ot\g^{1,0}$ can be written as \begin{equation}\bmu=\suml{i,j}
\mu_j^i\oom^j\ot T_i+\suml{i,\al}\mu_\al^i\oom^\al\ot T_i+\suml{j,\be}
\mu_j^\be\oom^j\ot W_\be+\suml{\al,\be}\mu_\al^\be \oom^\al \ot W_\be. 
\label{star}\end{equation} By Lemma~\ref{holomorphic forms} and
(\ref{eq:dbar}), \begin{eqnarray*} -\dbar\bmu &=& \suml{i,j}\mu_j^i\oom^j\we
\dbar T_i+\suml{i,\al}\mu_\al^i\oom^\al\we\dbar T_i\\&=& \suml{i,j,k,\be}
\mu_j^i E_{ki}^\be \oom^j \we \oom^k \ot W_\be+\suml{i,k,\al,\be}\mu_\al^i
E_{ki}^\be \oom^\al\we\oom^k \ot W_\be.\end{eqnarray*} This calculation gives
us a necessary and sufficient condition for $\bmu$ to be $\dbar$-closed, which
we now record as

\begin{lemma}\label{lemma:closed form} Suppose that an element $\bmu$ in
$\g^{*(0,1)}\ot\g^{1,0}$ is given by formula (\ref{star}). Then $\dbar \bmu=0$
if and only if \[\textstyle\suml i(\mu_j^i E_{ki}^\al-\mu_k^iE_{ji}^\al)=0
\qand \suml i\mu_\al^i E_{ji}^\be=0,\] for each $j,k,\al,\be$.\end{lemma}
\medbreak

We have a sequence \[ 0\to\g^{1,0}\to\g^{*(0,1)}\ot\g^{1,0}\to\ \cdots\ \to
\g^{*(0,k-1)}\ot\g^{1,0} \stackrel{\dbar_{k-1}}\to \g^{*(0,k)}\ot\g^{1,0}
\stackrel{\dbar_k}\to\cdots \] The next result comes as no surprise,
reflecting as it does the fact that our $\dbar$ operators are the natural ones
induced on invariant differential forms.

\begin{lemma}\label{seq} The above sequence is a complex, i.e.\ 
$\dbar_k\circ\dbar_{k-1}=0$ for all $k\ge1$.\end{lemma}

\bproof It suffices to verify the lemma for $k=1$. Let $\{\oom^p\},\{\oe_q\}$
be dual bases of $\g^{*(1,0)}$ and $\g^{1,0}$, where the indices $p,q$ run
over the entire range $1,\ldots,n+m$. By definition,\[\dbar V=\suml p\oom^p\ot
[\oe_p,V]^{1,0}.\] Applying $\dbar$ again,\[\dbar^2V=\suml p\dbar\oom^p\ot
[\oe_p,V]^{1,0}-\suml{p,q}(\oom^p\we\oom^q)\ot[\oe_q,[\oe_p,V]^{1,0}]^{1,0}.\]
Since $[\oe_p,\oe_q]^{1,0}=0$, we can delete the penultimate projection
$\>^{1,0}\>$ above. The Jacobi identity \[[\oe_p,[\oe_q,V]]-[\oe_q,[\oe_p,V]]
=[[\oe_p,\oe_q],V]\] implies that \[\suml{p,q}\oom^p\we\oom^q\ot
[\oe_q,[\oe_p,V]]^{1,0}=-\half\suml{p,q}\oom^p\we\oom^q\ot
[[\oe_p,\oe_q],V]^{1,0}.\]

If $\si$ is a $(1,0)$-form, we can now contract $\si$ with $\dbar^2V$ to
obtain the following form of type $(0,2)$: \[\begin{array}{rcl} \si(\dbar^2V)
&=& \suml p\si([\oe_p,V])\dbar\oom^p+\half\suml{p,q}\si([[\oe_p,\oe_q],V])
(\oom^p\we\oom^q)\\[10pt] &=&-2\suml p d\si(\oe_p,V)\dbar\oom^p-\suml{p,q}
d\si([\oe_p,\oe_q],V)(\oom^p\we\oom^q).\end{array}\] For this to vanish for
all $V$ and $\si$, we need to show that \[\textstyle 2\sum\limits_{r=1}^{n+m}
\dbar\oom^r\ot \oe_r=-\suml{p,q}(\oom^p\we\oom^q)\ot[\oe_p,\oe_q].\] This
equation amounts to stating that the $\oom^p\we\oom^q$ component of
$2\dbar\oe^r$ equals \[\oom^r([\oe_p,\oe_q])=2d\oom^r(\oe_p,\oe_q)=2\dbar
\oom^r(\oe_p,\oe_q),\] which is correct.\eproof

\begin{definition} Define $H^k_{\dbar}(\g^{1,0})$ to be the $k$th cohomology
$\ker\dbar_k/\image\dbar_{k-1}$ of the above complex; more precisely,
\[H^k_{\dbar}(\g^{1,0})=\frac{\ker\left(\dbar_k:\g^{*(0,k)}\ot\g^{1,0}\to
\g^{*(0,k+1)}\ot\g^{1,0}\right)}{\dbar_{k-1}\left(\g^{*(0,k-1)}\ot\g^{1,0}
\right)}.\] \end{definition}

We shall interpret these spaces geometrically in the next subsection.

\subsection{Dolbeault cohomology}

Let $\Ga$ be a $J$-invariant co-compact lattice in $G$, and $X=\Ga\bs
G$ the associated nilmanifold parameterizing left cosets. Let $\psi\colon G\to
G/C$ be the quotient map, where $C$ is the center of $G$. Since $G$ is 2-step
abelian, $G/C$ is abelian. In terms of the abelian varieties $F:=C/C\cap\Ga$
and $M:=\psi(G)/\psi(\Ga)$, we obtain a holomorphic fibration \[\Psi\colon
X\longrightarrow M\] with fiber $F$.

\begin{lemma} Let $\co_X$ and $\Th_X$ be the structure sheaf and the tangent
sheaf of $X$. For $p\ge1$, the direct image sheaves with respect to $\Psi$ are
\[\begin{array}{l} R^p\Psi_*\co_X=\We^p\fc^{*(0,1)}\ot\co_M=\fc^{*(0,p)}\ot
\co_M,\\[5pt] R^p\Psi_*\Psi^*\Th_M=\We^p\fc^{*(0,1)}\ot \Th_M
=\fc^{*(0,p)}\ot \Th_M.\end{array}\]\end{lemma}

\bproof The second identity is a consequence of the first, and the projection
formula. To prove the first, note that for any point $m$ in $M$,
\[(R^p\Psi_*\co_X)_m=H^p(\Psi^{-1}(m),\co_X)\cong H^p(C,\co_X).\] This has
constant rank and, by Grauert's Theorem, the direct image sheaf is locally
free. As $\Psi^{-1}(m)$ is isomorphic to a complex torus, for all
$p\ge1$,\[H^p(\Psi^{-1}(m),\co_X)=\We^pH^1(\Psi^{-1}(m),\co_X),\] The vector
bundle $R^p\Psi_*\co_X$ is isomorphic to $\We^pR^1\Psi_*\co_X$. Since the
space of vertical (0,1)-forms is trivialized by the left-invariant (0,1)-forms
given in Lemma~\ref{holomorphic forms}, we have \[R^p\Psi_*\co_X
\cong\We^pR^1\Psi_*\co_X\cong\We^p\fc^{*(0,1)}\ot\co_M,\] as required.\eproof

\begin{lemma}\label{structure sheaf} Let $\co_X$ and $\Th_X$ be the
structure sheaf and the tangent sheaf of $X$. Then \[\begin{array}{l} H^k(X,
\co_X)=\We^k\g^{*(0,1)}=\g^{*(0,k)},\\ H^k(X,\Psi^*\Th_M)=
\We^k(\g^{*(0,1)}) \ot \ft^{1,0}=\g^{*(0,k)} \ot \ft^{1,0}.\end{array}\]
\end{lemma}

\bproof Consider the Leray spectral sequence with respect to the
$\dbar$-operator and the holomorphic projection $\Psi$. One has
\[E_2^{p,q}=H^p(M, R^q\Psi_*\co_X),\qq E_{\infty}^{p, q}\Rightarrow
H^{p+q}(X,\co_X).\] From the previous lemma, when $q\ge1$, \begin{eqnarray*}
E_2^{p,q} &=&H^p(M,\We^q \fc^{*(0,1)} \ot\co_M)=\We^q \fc^{*(0,1)} \ot H^p(M,
\co_M)\\ &=&\We^q \fc^{*(0,1)} \ot H^p(M,\co_M)=\We^q \fc^{*(0,1)}\ot\We^p
\ft^{*(0,1)}.\end{eqnarray*} Note that every element in $E_2^{p,q}$ is a
linear combination of the tensor products of vertical $(0,q)$-forms and
$(0,p)$-forms lifted from the base. Since these forms are globally defined and
the differential $d_2$ is generated by the $\dbar$-operator, we have
$d_2=0$. It follows that the Leray spectral sequence degenerates at the
$E_2$-level. Therefore,\[H^k(X,\co_X)=\bigoplus_{p+q=k}E_2^{p, q}=\We^k
(\fc^{*(0,1)}\op \ft^{*(0,1)})=\We^k\g^{*(0,1)}.\]

Next, the spectral sequence for $\Psi^*\Th_M$ gives \[E_2^{p,q}=H^p(M,
R^q\Psi_*\Psi^*\Th_M),\qq E_{\infty}^{p, q}\Rightarrow H^{p+q}(X,\Psi^*\Th_M).
\] Moreover, $E_2^{p,q}$ is equal to \[H^p(M,\We^q\fc^{*(0,1)}
\ot\Th_M)=\We^q \fc^{*(0,1)}\ot H^p(M, \Th_M)=\We^q\fc^{*(0,1)}\ot\We^p
\ft^{*(0,1)}\ot \ft^{1,0}.\] Elements in $\ft^{1,0}$ are holomorphic vector
fields on $M$ and hence globally defined sections of $\Psi^*\Th_M$ on
$X$. Elements in $\We^q\fc^{*(0,1)}$ are pulled back to globally defined
$(0,q)$-forms on $X$. Crucially, elements in $\We^p\ft^{*(0,1)}$ are
globally-defined holomorphic $(0,p)$-forms on $X$, and the operator $d_2$ is
identically zero. Therefore, the spectral sequence degenerates at $E_2$. We
have \[H^k(X,\Psi^*\Th_M)=\bigoplus_{p+q=k}E_2^{p, q}=\We^k(\fc^{*(0,1)}\op
\ft^{*(0,1)})\ot \ft^{1,0}=\We^k(\g^{*(0,1)})\ot \ft^{1,0},\] as required.
\eproof

\begin{theorem}\label{identification} Let $X$ be a 2-step nilmanifold 
with an abelian complex structure. There is a natural isomorphism
$H^k(X,\Th_X) \cong H^k_{\dbar}(\g^{1,0})$.\end{theorem}

\bproof On the manifold $X$, we have the exact sequence \[ 0\to\fc^{1,0}\ot
\co_X \to \Th_X \to \Psi^*\Th_M\to 0.\] A piece of the corresponding long
exact sequence is \[\to \fc^{1,0}\ot H^k(X,\co_X)\to H^k(X,\Th_X) \to
H^k(X,\Psi^*\Th_M) \stackrel{\delta_k}\to\fc^{1,0}\ot H^{k+1}(X,\co_X) \to\]
>From the last section, the coboundary map is \[\delta_k:\g^{*(0,k)}\ot
\ft^{1,0}\to\g^{*(0,k+1)}\ot \fc^{1,0},\] and so \[H^k(X,\Th_X)\cong\ker
\delta_k\op\frac{\g^{*(0,k)}\ot\fc^{1,0}}{\delta_{k-1}(\g^{*(0,k-1)} \ot
\ft^{1,0})}.\]\vs

We calculate the coboundary maps by chasing the commutative diagram\vs

\[\begin{array}[c]{ccccccccc} 0 & \to & T^{\ast(0,k+1)}\ot\fc^{1,0} & \to & 
T^{\ast(0,k+1)}\ot\Th_X & \to & T^{\ast(0,k+1)}\ot \Psi^{\ast}\Th_{M} & \to &
0\\[5pt] & & \big\uparrow & & \big\uparrow & & \big\uparrow & & \\[8pt] 0 &
\to & T^{\ast(0,k)}\ot \fc^{1,0} & \to & T^{\ast (0,k)}\ot\Th_{X} & \to &
T^{\ast(0,k)}\ot\Psi^{\ast} \Th_{M} & \to & 0\end{array}\]\\ The vertical maps
are $\dbar$'s for Dolbeault cohomology. More specifically, if $\nabla$ is the
Chern connection, if $\om$ is a $(0,k)$-form and $V$ a vector field of type
$(1,0)$, then \[\dbar(\om \ot V)=\dbar\om\ot V+(-1)^k\om\we\dbar^\nabla V.\]
In the following computation, we let $\{e_p\}$ be a left-invariant basis for
$\g^{1,0}$ and $\{\om^p\}$ the dual basis.

Let $\oom$ be a $(0,k)$-form. Let $V$ be an element in $\ft^{1,0}$, considered
as a holomorphic vector field on $M$ and a holomorphic section of $\Psi^*
\Th_M$. Let $\tilde{V}$ a smooth lifting of this section to be a section of
$\Th_X$. Then \[\delta_k(\oom\ot V)=\dbar \oom\ot\tilde V+ (-1)^k\suml p\oom
\we\oom^p\ot[\ol e_p,\tilde V]^{1,0}.\] The element $T_j$ in $\g_\bC$ could be
considered as holomorphic vector field on $M$. It could also be considered as
a smooth vector field on $X$. Considering the latter a lifting of the former
and applying the above formula, we see that $\delta_k=\dbar_k$ on $\g^{*(0,k)}
\ot\ft^{1,0}$. Now $\dbar_k (\g^{*(0,k)} \ot \fc^{1,0})=0$, and so \[
\delta_{k-1}(\g^{*(0,k-1)}\ot \ft^{1,0})=\dbar_{k-1}(\g^{*(0,k-1)}\ot
\g^{1,0}).\] Since the Lie algebra $\g$ is 2-step nilpotent, $\image
\dbar_{k-1}\subseteq\g^{*(0,k)}\ot \fc^{1,0}$. Also, we have \[\ker\delta_k=
\ker\dbar_k\cap (\g^{*(0,k)} \ot \ft^{1,0}).\] Therefore, \begin{eqnarray}
H^k(X,\Th_X) &=& \ker\dbar_k\cap(\g^{*(0,k)}\ot\ft^{1,0})\op\frac{\g^{*(0,k)}
\ot\fc^{1,0}}{\image\dbar_{k-1}} \nonumber\\[8pt] &=& \frac{ \ker\dbar_k\cap
(\g^{*(0,k)}\ot\ft^{1,0}\op\g^{*(0,k)}\ot\fc^{1,0})}{\image\dbar_{k-1}}=
H^k_{\dbar}(\g),\end{eqnarray} as stated.\eproof

To summarize the results in this section, we shall say that a tensor on $X$ is
{\em invariant} if its pull-back to $G$ by the quotient map is invariant by
left-translation by $G$. Lemma \ref{structure sheaf} and Theorem
\ref{identification} then allow us to formulate

\begin{theorem}\label{invariant element} The Dolbeault cohomology on $X$ with
coefficients in the structure and tangent sheaf can be computed using
invariant forms and invariant vectors.\end{theorem}

Although the above proof relies on the 2-step property, one might expect that
this result has a more generally validity, at least in the nilpotent context.
For an independent approach to this problem, see \cite{CF}.

\section{Deformation theory}

We shall shortly be in a position to apply the Kuranishi method to construct
deformations. But first, we shall exhibit harmonic representatives in the
Dolbeault cohomology groups.

\subsection{Harmonic theory}

Theorem \ref{identification} reduces the question to finite-dimensional vector
spaces, and we may choose an invariant Hermitian structure on $X$ of the type
mentioned mentioned after Lemma~\ref{holomorphic forms}. We use the resulting
inner product on $\g^{*(0,k)}\ot\g^{1,0}$ to define the orthogonal complement
of $\image \dbar_{k-1}$ in $\ker\dbar_k$. Denote this space by
$\image^\perp\dbar_{k-1}$.

\begin{theorem}\label{harmonic rep} The space $\image^\perp\dbar_{k-1}$ is
a space of harmonic representatives for the Dolbeault cohomology
$H^k(X,\Th_X)$ on the compact complex manifold $X$.\end{theorem}

\bproof It suffices to prove that an element \begin{equation}\label{imperp}
\suml p\osi^p\ot e_p\in\image^\perp\dbar_{k-1}\subseteq\g^{*(0,k)}\ot\g^{1,0}
\end{equation} is $\dbar^*$-closed on the manifold $X$.

Any section of the trivial bundle over $X$ with fibre $\g^{*(0,k-1)}\ot
\g^{1,0}$ is a sum of elements of the type $f\oet\ot V$, where $f$ is a smooth
function, $\oet\in\g^{*(0,k-1)}$ and $V\in\g^{(1,0)}$. By
Lemma~\ref{holomorphic forms}, $\osi^p$ and $\oet$ are $\dbar$-closed. Using
double angular brackets for the $L^2$ inner product and summing over repeated
indices, we calculate \begin{eqnarray*} \ll\dbar^*(\osi^p\ot e_p),\ 
f\oet\ot V\rr &=& \ll\osi^p,\dbar(f\oet)\rr\langle e_p,V\rangle+
(-1)^{k-1}\ll\osi^p\ot e_p,\ f\oet\we\dbar V\rr\\[3pt]&=& \ll\dbar^*\osi^p,
f\oet\rr\langle e_p,V\rangle+(-1)^{k-1}\ll\osi^p\ot e_p,\ f\oet\we
\dbar V\rr.\end{eqnarray*}

The basis $\{\om^i,\om^\al\}$ of Lemma~\ref{holomorphic forms} determines
a complex volume form that we may use to identify $\dbar^*$ with
$\pm\ast\!\dbar\,\ast$, where $\ast$ is the corresponding $SU(n+m)$ invariant
antilinear mapping $\g^{*(0,k)}\to \g^{*(0,n+m-k)}$. It follows that
$\dbar^*\osi=0$. The remaining term \[\ll\osi^p\ot e_p,\ f\oet\we\dbar
V\rr= \int_X \ol f\langle\osi^p\ot e_p,\oet\we\dbar V\rangle=
\langle\osi^p\ot e_p, \dbar(\oet\ot V)\rangle\int_X \ol f\] vanishes by
assumption (\ref{imperp}).\eproof

\begin{corollary}\label{dbar-*} Let $\bmu\in\g^{*(0,k)}\ot\g^{1,0}$. Then
$\dbar^*\bmu$ with respect to the $L_2$-norm on the compact manifold $X$ is
equal to $\dbar^*\bmu$ with respect to the Hermitian inner product on the
finite-dimensional vector spaces $\g^{*(0,k)}\ot\g^{1,0}$.\end{corollary}

\bproof This follows from the displayed formulae in the previous proof.
\eproof

\subsection{The Schouten-Nijenhuis bracket}

If $\oom\ot V$ and $\oom'\ot V'$ are vector-valued (0,1)-forms representing
elements in $H^1(X,\Th_X)$, their product with respect to the
Schouten-Nijenhuis bracket is a vector-valued (0,2)-form \[\{\cdot,\cdot\}:
H^1(X,\Th_X)\times H^1(X,\Th_X)\to H^2(X,\Th_X).\] It is defined at the level
of forms by \[\{\oom\ot V,\oom'\ot V'\}=\oom'\we L_{V'}\oom\ot V +\oom\we
L_V\oom'\ot V' +\oom\we\oom' \ot [V, V'].\] Via the isomorphism with Lie
algebra cohomology, elements in $H^1(X,\Th_X)$ lie in $\image^\perp\dbar_0$.
Since the vector and form parts are all left-invariant, $\iota_V\oom'$ is a
constant. Therefore, $L_V\oom'=d\iota_V\oom'+\iota_V d\oom'=\iota_V d\oom'$,
and \[\{\oom\ot V,\oom'\ot V'\}=\oom'\we \iota_{V'}d\oom\ot V +\oom\we \iota_V
d\oom'\ot V' +\oom\we\oom' \ot [V, V'].\] The complex structure is abelian, so
$[V, V']=0$ for all $(1,0)$-vectors, and \begin{equation}\label{Schouten}
\{\oom\ot V, \oom'\ot V'\}=\oom'\we \iota_{V'}d\oom\ot V +\oom\we\iota_V
d\oom'\ot V'.\end{equation}

Using the vector space direct sum $\g=\ft\op\fc$, we write
\begin{equation}\label{4 terms} \g^{*(0,1)}\ot\g^{1,0}=(\ft^{*(0,1)}\ot
\ft^{1,0})\op (\fc^{*(0,1)}\ot \ft^{1,0})\op(\fc^{*(0,1)}\ot \fc^{1,0})\op
(\ft^{*(0,1)}\ot \fc^{1,0}).\end{equation} If $\oom\ot V\in\ft^{*(0,1)}
\ot\fc^{1,0}$ then $d\oom=0$, because all elements in $\ft^{*(k,l)}$ are
closed. On the other hand, $d\oom'\in\ft^{*(1,1)}$. Since $\iota_V d\oom'=0$
for $V\in\fc^{1,0}$, we have \begin{equation}\label{first zero}
\{\ft^{*(0,1)}\ot \fc^{1,0},\ \g^{*(0,1)}\ot \g^{1,0}\}=0. \end{equation}

In order to compute $\{\bmu,\bnu\}$ on $\image^\perp\dbar_0$, we compute the
bracket amongst elements in the obvious basis. In view of (\ref{first zero}),
we need to calculate the brackets arising from the first three summands in
(\ref{4 terms}). There are six types of bracket to calculate. Since $\oom^k$
and $\oom^j$ are closed, \[\{\oom^j\ot T_i,\oom^k\ot T_l\}=0\] 
\[\{\oom^j\ot T_i,\oom^\al\ot T_l\}=\oom^j\we\iota_{T_i}d\oom^\al\ot T_l=
-\ol E_{ik}^\al\oom^j\we\oom^k\ot T_l\] \begin{equation}\label{eq:basic terms} 
\{\oom^j\ot T_i,\oom^\al\ot W_\si\}=\oom^j\we\iota_{T_i}d\oom^\al\ot W_\si=
-{\ol E}_{ih}^\al\oom^j\we\oom^h\ot W_\si \end{equation} \[\{\oom^\al\ot T_l,
\oom^\be\ot T_j\}=-\ol E_{lh}^\be\oom^\al\we\oom^h\ot T_j-\ol E_{jh}^\al
\oom^\be\we\oom^h\ot T_l \] \[\{\oom^\al\ot T_l,\oom^\be\ot W_\ga\}=\oom^\al\we
\iota_{T_l}d\oom^\be\ot W_\ga=-\ol E_{lh}^\be\oom^\al\we\oom^h\ot W_\ga \]
\[\{\oom^\al\ot W_\be,\oom^\ga\ot W_\delta\}=0.\]

The above formulae allow us to calculate $\{\bmu,\bnu\}$. If $\bmu$ is given in
coordinates as in (\ref{star}) and $\bnu$ similarly then, suppressing summation
signs, we have \begin{eqnarray} \{\bmu,\bnu\}&=&-(\mu_j^i\nu_\al^\ell+\nu_j^i
\mu_\al^\ell)\ol E_{ik}^\al\,\oom^j\we\oom^k\ot T_\ell\nonumber\\&&
-(\mu_\al^\ell\nu_\be^j+\nu_\al^\ell\mu_\be^j)(\ol E_{\ell k}^\be
\oom^\al\we\oom^k\ot T_j +\ol E_{jk}^\al\oom^\be \we\oom^k\ot T_\ell)
\nonumber\\&& -(\mu_j^i\nu_\ga^\delta+\nu_j^i\mu_\ga^\delta)\ol E_{ik}^\ga\,
\oom^j\we \oom^k\ot W_\delta\\
&&-(\mu_\al^i\nu_\ga^\delta+\nu_\al^i\mu_\ga^\delta)\ol
E_{ik}^\ga\,\oom^\al\we\oom^k\ot W_\delta\label{eq:generic bracket}
\end{eqnarray} In particular, \begin{eqnarray}\{\bmu,\bmu\} &\!\!=\!\!&
-2\mu_j^i\mu_\al^\ell\ol E_{ik}^\al\,\oom^j\we\oom^k\ot T_\ell\nonumber\\&&
-2\mu_\al^\ell\mu_\be^j (\ol E_{\ell k}^\be\oom^\al\we\oom^k\ot T_j+\ol
E_{jk}^\al\oom^\be\we\oom^k\ot T_\ell )\\&& -2\mu_j^i\mu_\ga^\delta\ol
E_{ik}^\ga\,\oom^j\we\oom^k\ot W_\delta -2\mu_\al^i\mu_\ga^\delta\ol
E_{ik}^\ga\,\oom^\al\we\oom^k\ot W_\delta,\nonumber\label{eq:self}
\end{eqnarray} which is of course is an element of
$\g^{*(0,2)} \ot\g^{1,0}$.

\subsection{Kuranishi theory}

To construct deformations, we apply Kuranishi's recursive formula. Let
$\{\be_1,\dots,\be_N\}$ be an orthonormal basis of the harmonic
representatives of $H^1(X,\Th_X)$. For any vector $\bt=(t_1,\dots, t_N)$
in $\bC^N$, let \begin{equation}\label{mut}\bmu(\bt)=t_1 \be_1+\cdots+t_N\be_N.
\end{equation} We set $\bph_1=\bmu$, and next define $\bph_r$ inductively for 
$r\ge2$.

Consider the $\dbar$-operator on $X$ with respect to the Hermitian metric $h$
previously defined, its adjoint operator $\dbar^*$, and the Laplacian
\begin{equation}\label{tri}\triangle=\dbar\dbar^*+\dbar^*\dbar.\end{equation}
Let $\cg$ be the corresponding Green's operator that inverts $\triangle$ on
the orthogonal complement of the space of harmonic forms, and let $\{\ ,\ \}$
denote the Schouten-Nijenhuis bracket. Then we set \begin{equation}\label{phi}
\bph_r(\bt)=\half\sum_{s=1}^{r-1}\dbar^*\cg\{\bph_s(\bt),\bph_{rs}(\bt)\}=
\half\sum_{s=1}^{r-1}\cg\dbar^*\{\bph_s(\bt),\bph_{r-s}(\bt)\},\end{equation}
and consider the formal sum \begin{equation}\label{fsum}\bPh(\bt)=
\sum_{r\ge1}\bph_r.\end{equation}

Let $\{\ga_1,\dots,\ga_M\}$ be an orthonormal basis for the space of harmonic
$(0, 2)$-forms with values in $\Th_X$. Define $f_k(\bt)$ to be the $L^2$-inner
product $\ll\{\bPh(\bt),\bPh(\bt)\},\ga_k\rr$. Kuranishi theory asserts the
existence of $\eps>0$ such that \begin{equation}\label{family}
\{\bt\in\bC^N:|\bt|<\eps,\ f_1(\bt)=0,\dots,f_M(\bt)=0\}\end{equation} forms a
locally complete family of deformations of $X$. We shall denote this set by
$\Kur$. For each $\bt\in\Kur$, the associated sum $\bPh=\bPh(\bt)$ satisfies
the integrability condition \begin{equation}\label{MC}\dbar\bPh+\half\{\bPh,
\bPh\}=0\end{equation} that now follows from (\ref{phi}) and the definition of
$\cg$.

More explicitly, we may treat $\bPh$ is a linear map from $(0,1)$-vectors to
$(1,0)$-vectors. It determines a complex structure on our manifold $X$ whose
distribution of (0,1)-vectors is given by \begin{equation}\label{Abel}
\left\{\begin{array}{l} \ol S_j=\ol T_j+\bPh(\ol T_j),\\[3pt]\ol V_\al=\ol
W_\al+\bPh(\ol W_\al).\end{array}\right.\end{equation} This set of equations
is analogous to the gauge-theoretic defininition of a connection as $d_A=d+A$,
where $A$ is a matrix of 1-forms. In principal bundle language, $d_A$
determines a horizontal distribution formed from the flat one by adding A as a
vertical component. Then (\ref{MC}) is the analogue of setting the curvature
of $d_A$ to be zero, and assures us that the new distribution (\ref{Abel}) is
closed under Lie bracket.

We are now ready to make precise the first theorem of the Introduction:

\begin{theorem}\label{main} Let $G$ be a 2-step nilpotent Lie group with
co-compact subgroup $\Ga$, and let $J$ be an abelian invariant complex
structure on $X=\Ga\bs G$. Then the deformations arising from $J$
parameterized by (\ref{family}) are all invariant complex structures.
\end{theorem}

\bproof It suffices to show that every term in the power series (\ref{fsum})
lies in $\g^{*(0,1)}\ot \g^{1,0}$. We shall prove this by induction. By
Theorem \ref{invariant element}, $\bph_1=\bmu$ belongs to this space.

Assume that $\bph_s\in\g^{*(0,1)}\ot \g^{1,0}$ for all $1\le s\le r-1$. The
computations of \S4.2 show that $\{\bph_s,\bph_{r-s}\}$ is always contained in
\[\g^{*(0,2)}\ot\g^{1,0}=\g^{*(0,2)}\ot\g^{1,0}=\image\dbar_1\op
\image^\perp\dbar_1.\] Let $\pi_0$ denote projection to the subspace
$\image\dbar_1$.

By Theorem \ref{harmonic rep}, the component $\image^\perp\dbar_1$ is the
harmonic part of $H^2(X,\Th)$. Since (\ref{tri}) satisfies $\triangle\circ
\dbar=\dbar\circ\triangle$, and $\image\dbar_1$ is orthogonal to the harmonic
part, $\triangle$ maps $\image\dbar_1$ isomorphically onto itself. It follows
that \[\cg\{\bph_s,\bph_{r-s}\}=\cg\pi_0\{\bph_s,\bph_{r-s}\}\subseteq
\image\dbar_1.\] In particular, it is an invariant tensor. Corollary
\ref{dbar-*} shows that $\dbar^*{\cg}\{\bph_s,\bph_{r-s}\}$ is again an
invariant tensor. The same is true of $\bph_r$. By induction, (\ref{fsum}) is
an infinite series of invariant tensors. \eproof

\section{Deformations leading to abelian structures}

In the light of Theorem \ref{main}, we are now ready to identify deformations
of $J$ leading not just to invariant complex structures, but to \emph{abelian}
ones.

Given an element $\bmu=\bmu(\bt)$ in the virtual parameter space
$H^1(X,\Th_X)$ as in (\ref{mut}), we apply the preceding method to generate
the infinite series (\ref{fsum}), and consider (\ref{Abel}). For the latter to
define an \emph{abelian} complex structure, the Lie bracket of any pair of
$(0,1)$-vectors must in fact vanish identically. In this case, we shall say
that $\bmu$ generates an \emph{abelian deformation}.  Such an assumption leads
to the following equations: \begin{eqnarray} \label{eq:lie1} \lb\ol S_j,\ol
S_k \rb &=& 0,\qq 1 \le j,k \le n \\\label{eq:lie2} \lb\ol S_j,\ol V_\al\rb
&=& 0,\qq 1\le j\le n,\q n+1\le \al \le n+m\\\label{eq:lie3}\lb\ol V_\al,\ol
V_\be \rb &=& 0,\qq n+1\le\al,\be \le n+m.\end{eqnarray} Since $\ol W_\al$ is
in the center, \[[\ol V_\al,\ol V_\be]= [\bPh(\ol W_\al),\bPh(\ol W_\be)],\]
and this vanishes since the original complex structure is abelian. Therefore,
equation (\ref{eq:lie3}) is satisfied automatically.

Let us examine the infinitesimal consequence of the first two equations. Let
$t$ represent a real variable, and replace $\bmu$ by $t\bmu$ so that $\bPh$
becomes $\sum t^r\bph_r$. Then in the notation of (\ref{star}), equation
(\ref{eq:lie1}) leads to \begin{eqnarray} 0 &=& \frac d{dt}\Big|_{t=0}\lb\ol
S_j,\ol S_k\rb=\lb\ol T_j,\bph_1\ol T_k\rb + \lb\bph_1\ol T_j,\ol T_k\rb=
\lb\ol T_j,\mu_k^iT_i\rb+\lb \mu_j^iT_i,\ol T_k\rb\nonumber\\&=&
(\mu_k^iE_{ji}^\al-\mu_j^iE_{ki}^\al)W_\al+(\mu_k^iF_{ji}^\al-\mu_j^i
F_{ki}^\al)\ol W_\al.\label{eq:infin1}\end{eqnarray} The coefficient of
$W_\al$ vanishes when $\dbar\bmu=0$, by Lemma~\ref{lemma:closed form}.
Equation (\ref{eq:lie2}) leads to \begin{eqnarray} 0 &=& \frac d{dt}
\Big|_{t=0} \lb\ol S_j,\ol V_\al\rb=\lb\ol T_j,\bph_1\ol W_\al\rb=\lb\ol
T_j,\mu_\al^iT_i\rb \nonumber\\ &=& \mu_\al^iE_{ji}^\be W_\be+\mu_\al^i
F_{ji}^\be \ol W_\be.\label{eq:infin2} \end{eqnarray} The coefficient of
$W_\be$ is again $0$ when $\bmu$ is $\dbar$-closed.

The above calculations give a set of necessary conditions limiting the type
of deformations that one needs to consider. They motivate 

\begin{definition}\label{def:conditiona} A form $\bmu$ given in coordinates as 
in (\ref{star}) satisfies {\em Condition A} if \[\textstyle \suml i(\mu_j^i
F_{ki}^\al-\mu_k^iF_{ji}^\al)=0\qand\suml i\mu_\al^iF_{ji}^\be=0,\] for
each $j,k,\al,\be$.\end{definition}

It is striking that these conditions are completely analogous to those of
Lemma~\ref{lemma:closed form}. In view of (\ref{eq:infin1}) and
(\ref{eq:infin2}), we can now state

\begin{proposition}\label{prop:conditiona} A parameter $\bmu$ represents an 
infinitesimal abelian deformation if and only if it is $\dbar$-closed and
satisfies Condition A.\end{proposition}

Next suppose that $\bmu$ and $\bnu$ are vector-valued 1-forms that are both
$\dbar$-closed and satisfy Condition A. Since $\ol{E}_{ij}^\al
=-F_{ji}^\al$, every term in (\ref{eq:generic bracket}) is equal to $0$.
For example, the first term $-\mu_j^i\nu_\al^\ell \ol{E}_{ik}^\al
\oom^j\we \oom^k \ot T_\ell$ is equal to \[-\mu_j^i\nu_\al^\ell
F_{ki}^\al \oom^j\we\oom^k \ot T_\ell=-\nu_\al^\ell\left(\mu_j^i
F_{ki}^\al - \mu_k^iF_{ji}^\al\right) \oom^j \ot \oom^k \ot T_\ell=0.\]
and similarly every term in $\{\bmu,\bnu\}$ is equal to zero. In particular,
$\{\bmu,\bmu\}=0$.

Using the recursive formula (\ref{phi}), the higher order terms are all
equal to zero, and so the series $\bPh$ and $\bmu$ coincide by
construction. Furthermore, $\{\bPh,\bPh\}=\{\bmu,\bmu\}=0$, and there is no
additional obstruction to integrability. Therefore,

\begin{proposition}\label{integrable} On a 2-step nilmanifold $X$ with abelian
complex structure, an element in $H^1(X,\Th_X)$ is infinitesimally abelian
only if it is integrable to a 1-parameter family of abelian complex structures.
\end{proposition}

Our main result concerning the deformation of abelian complex structures is

\begin{theorem}\label{mainAb} On a 2-step nilmanifold with abelian complex
structure, a parameter $\bmu$ in $\g^{*(0,1)}\ot\g^{1,0}$ generates an abelian
deformation if and only if it is $\dbar$-closed and satisfies Condition~A. 
\end{theorem}
  
\bproof If $\bmu$ generates an abelian deformation, it is infinitesimally
abelian. By Proposition~\ref{prop:conditiona}, the form is $\dbar$-closed and
satisfies Condition A. 

Conversely, if $\bPh$ is $\dbar$-closed, it represents a cohomology class in
$H^1(X,\Th)$. Since it also satisfies Condition~A, it is infinitesimally
abelian. By Proposition \ref{integrable}, it represents an integrable abelian
complex structure.\eproof

\subsection{Fully abelian deformations}

We are curious to know when the entire virtual parameter space
$H^1(X,\Th_X)$ integrates to abelian complex structures.

\begin{theorem}\label{fully} Let $X=\Ga\bs G$ be a compact 2-step
nilmanifold endowed with an abelian complex structure. Suppose that every
direction of the virtual parameter space is integrable to a 1-parameter family
of abelian complex structures and that the dimension of the center of Lie
algebra $\g$ is equal to $1$. Then $\g$ is isomorphic to the direct sum of a
Heisenberg algebra and a 1-dimensional abelian algebra.\end{theorem}

\bproof Given the hypothesis on the center, we may as well drop the index
$\al$ in $E^\al_{ij}$. This feature makes the subsequent construction
possible.

Given the structural constants, for each set of $j,k,l,m$, choose an element
$\bmu$ in $\g^{*(0,1)}\ot \g^{1,0}$ by setting $\mu_k^l=E_{km}$ and $\mu_j^m=
E_{jl}$ and all other terms are set to zero. By Lemma \ref{lemma:closed form},
each such $\bmu$ is closed. By equation (\ref{eq:self}), such $\bmu$ satisfies
the equation $\{\bmu,\bmu\}=0$ and therefore there is no obstruction for it to
represent an integrable complex structure.

By hypothesis, $\bmu$ represents an abelian complex structure. By Condition A,
\[E_{km}\ol{E}_{lj}-E_{jl}\ol{E}_{mk}=0 .\] It follows that $\left|E_{km} 
\right|^2=\left| E_{mk}\right|^2$ for all $k$ and $m$.

If every $E_{km}$ vanishes then the algebra is abelian. On the other hand, if
at least one $E_{km}$ is non-zero, then $E_{mk}\ne0$. For every $E_{jl}\ne0$,
the ratio \[\frac{\ol E_{jl}}{E_{lj}}=\frac{\ol E_{km}}{E_{mk}}\] is
independent of the choice of $j,l$. Hence, there exists a real number $\theta$
such that \begin{equation}\label{e-symmetry} e^{i\theta}E_{jl}=\ol E_{lj},
\end{equation} for every pair of indices $(j,l)$. It follows that 
\[\lb\ol T_j,T_l\rb=E_{jl}W+F_{jl}\ol W=E_{jl}W-\ol E_{lj}\ol W=E_{jl}W-
e^{i\theta }E_{jl}\ol{W}.\] Choosing \[ D_{jl}=e^{i(\pi+\theta)/2}E_{jl},\q
U=e^{-i(\pi+\theta)/2}W\] gives \[\lb\ol T_j,T_l\rb=D_{jl}(U+\ol U).\] With
(\ref{e-symmetry}), we find that the matrix $(D_{jl})$ is skew-Hermitian. If
we now choose a basis of $\left( 0,1\right) $-vectors so that the matrix $D$
is diagonal, the diagonal entries are purely imaginary or zero. The
restriction on the central dimension forces the matrix $D$ to be a constant
multiple of the identity matrix. It follows that the structural equations
exactly mirror those of the Heisenberg algebra as seen in
(\ref{hei-structure}).\eproof

\begin{exa} There exist examples satisfying the first hypothesis of the
theorem, but not the second. To see this, take $\g$ to be the real
8-dimensional Lie algebra with non-zero complex structural equations \[\lb\ol
T_1,T_1\rb=W_3+\ol W_3, \q\lb\ol T_2,T_2\rb=W_4 +\ol W_4,\] and real
4-dimensional center. By Lemma \ref{lemma:closed form}, $\mu_1^2=\mu_2^1=0$,
and $\mu_\al^h=0$ for $1\le h\le 2$ and $3\le\al\le4$. It follows that
\[H^1(X,\Th_X)=\langle \oom^1\ot T_1,\> \oom^2 \ot T_2,\>\oom^1 \ot W_4,\>
\oom^2\ot W_3 \rangle,\] and one may check that each direction is integrable
to abelian complex structures. Globally, the associated compact complex
manifold is the product of two primary Kodaira surfaces.\end{exa}

\section{Six-dimensional structures}

In dimension 6, there are precisely six classes of 2-step groups or
nilmanifolds with an abelian complex structure \cite{simon}. Namely, the
abelian group $R^6$, the product $H_5\times R^1$ of a 5-dimensional Heisenberg
group with a 1-dimensional group, the product $H_3\times R^3$ of the
3-dimensional Heisenberg group with a 3-dimensional abelian group, the product
$H_3\times H_3$ of two 3-dimensional Heisenberg groups, and the groups $P_6$
and $W_6$. These were encountered in \S2.1.

We shall use Example~\ref{W6} to illustrate that results in this article
produce adequate information for finding the parameters for integrable and
abelian deformations. Using (\ref{J1}), consider the basis \[T_1=X_1-iX_2,\q
T_2=X_3+iX_4,\q W=Z_5+iZ_6\] of $\g^{1,0}$; the corresponding basis for
$\g^{*(0,1)}$ is $\{\oom^1,\oom^2,\oom\}$. The structural equations yield \[
[\ol T_1,T_2]=-W,\] so that \[ E_{12}=-1,\q F_{21}=1,\] and all other
structural constants are equal to zero. In particular,
\begin{equation}\label{ii} d\oom=\omega^1\we\oom^2,\q
\iota_{T_1}d\oom=\oom^2,\q\iota_{T_2}d\oom=0.\end{equation}

Mimicking the proof of Lemma~\ref{lemma:closed form}, any element $\bmu\in
\g^{*(0,1)}\ot\g^{1,0}$ can be written as \[\bPh=\mu_j^i\oom^j \ot T_i+ 
\mu^i_3\oom\ot T_i+\mu_j^3\oom^j\ot W+ \mu_3^3\oom \ot W,\] and
\begin{equation}\label{dbar-phi}\dbar\bmu=\mu_2^2\oom^2\we\oom^1\ot W
+\mu^2_3\oom\we\oom^1 \ot W.\end{equation} This shows that $\bmu$ is closed if
and only if $\mu^2_2=\mu^2_3=0$. Since $\dbar T_2=-\oom^1\ot W$, the space of
harmonic elements is in the orthogonal complement of $\oom^1\ot W$. Therefore,
\[ \dim H^1(X,\Th_X)=\dim\{\mu^2_2=\mu^2_3=\mu_1^3=0\}=6.\] Using
Definition~\ref{def:conditiona}, we see immediately that $\bmu$ satisfies
Condition A if and only if \begin{equation}\label{CondA}\mu^1_1=\mu^1_3=0.
\end{equation} The number of parameters corresponding to abelian deformations
is therefore 4.

Alternatively, we can count the number of integrable parameters, disregarding
the abelian issue. To do so, we employ the recursive formula from \S4.3, and
first calculate the self-bracket of a harmonic representative $\bmu$. Using
(\ref{ii}),(\ref{Schouten}), (\ref{eq:basic terms}), we deduce that
\[\begin{array}{ll}\{\bmu,\bmu\}\kern-8pt&=\{\mu^1_1\oom^1\ot T_1+\mu^1_3
\oom\ot T_1+\mu_3^3\oom\ot W,\ \mu^1_1\oom^1\ot T_1+\mu^1_3\oom\ot T_1+\mu_3^3
\oom\ot W\}\\[3pt] &= \mu^1_3(2\mu^1_1\oom^1\we\oom^2\ot T_1\!+\!\mu^1_3\oom
\we\oom^2\ot T_1\!+\!2\mu_3^3\oom\we\oom^2\ot W)\!-\!2\mu^1_1\mu_3^3\dbar(
\oom^2\ot T_2).\end{array}\] Using (\ref{phi}), we take \[\bph_2=\mu^1_1
\mu_3^3\oom^2\ot T_2.\] This quadratic correction term exactly corresponds to
the equation $d=-av$ in \cite[Proposition 4.2]{KS}.

If we set $\bPh=\bmu+\bph_2$, then (\ref{MC}) becomes \[
\mu^1_3(2\mu^1_1\oom^1\we\oom^2\ot T_1+\mu^1_3\oom\we\oom^2\ot T_1
+2\mu_3^3\oom\we\oom^2\ot W)=0.\] The resulting deformation is therefore
integrable if and only if $\mu^1_3=0$, so there is a total of 5 integrable
parameters. As predicted by Thereom~\ref{mainAb}, the obstruction $\mu^1_3$
already features in the abelian equations (\ref{CondA}).

\def\y{\\[5pt]}\def\FIG{
\vspace{15pt}\[\begin{array}{|c|c|c|c|c|c|}\cline{2-6} 
\multicolumn{1}{l|}{}& \vphantom{\int^{l^a}}
                      d & h^0 & h^1 &\dim\Kur &\dim\Abel\y\hline
T^6 \vphantom{\int^{l^a}} & \q 9\q  & \q 3\q  & \q 9\q & 9 & 9 \y
(\Ga\bs H_5) \times S^1   & 6 & 1 & 4 & 4 & 4 \y
(\Ga\bs H_3) \times T^3   & 7 & 2 & 6 & 6 & 6 \y
(\Ga\bs H_3)\!\times\!(\Ga\bs H_3) & 6 & 1 & 4 & 4 & 3 \y
\Ga\bs W_6                  & 6 & 2 & 6 & 5 & 4 \y
\Ga\bs P_6                  & 6 & 1 & 4 & 4 & 3 \y\hline
\end{array}\]\par
\centerline{Table}\vspace{15pt}}

Let $\g$ denote a real 6-dimensional nilpotent Lie algebra admitting a complex
structure. The table in \cite[Appendix]{simon} displays, for each such $\g$,
the complex dimension of the space $\cC(\g)$ of {\em invariant} complex
structures at a smooth point of one of its connected component. This was done
with little regard for when complex structures are equivalent, in the
knowledge that subsequent work would clarify the findings. The following table
compares these computations with results yielded by the techniques of this
paper.

The last five columns display the complex dimension

\begin{enumerate}\setlength{\itemsep}{-3pt}

\item $d$ of $\cC(\g)$, 

\item $h^0$ of the space $\dim H^0(X,\Th_X)$ of infinitesimal
automorphisms,

\item $h^1$ of the virtual parameter space $H^1(X,\Th_X)$,

\item of the space $\Kur$ of (\ref{family}), or the number of integrable
parameters,

\item of the subspace $\Abel$ of $\Kur$ describing abelian deformations.
\end{enumerate}

\n relative to the complex structures defined by (\ref{J}),(\ref{JJ}),(\ref
{JJJ}),(\ref{J1}),(\ref{JP}), for each of the last five rows in turn.

\FIG 

If $J$ is an invariant complex structure with unobstructed deformations on a
nilmanifold, $\cC(\g)$ has the same dimension as the kernel of \[\dbar:
\g^{*(0,1)}\ot\g^{1,0}\to \g^{*(0,2)}\ot\g^{1,0},\] whereas $\dim\Kur=h^1$.
Since the dimension of the image $\dbar(\g^{*(0,0)}\ot\g^{1,0})$ equals
$3-h^0$, we deduce further that \[d=3-h^0+h^1\] if $J$ is a generic point of
$\cC(\g)$.

At points of $\cC(\g)$ where $h^0$ jumps to a higher value, the Kuranishi
method is unable to detect the additional equivalences that come into play at
neighbouring points where the symmetry group drops. Consequently, we can only
assert that $d+h^0-3$ is an {\em upper bound} for $\dim\Kur$. In the Table,
these two numbers only disagree for $W_6$, and this is because $J_1$ was `too'
special a choice at which to carry out the computations. If we work instead at
a nearby point $J'$ corresponding to $\mu_3^3\ne0$, then $h^0=1$ and the
dimensions of $\Kur$ and $\Abel$ drop to 4 and 3. This is because the orbit
$J'\cdot W_6$ under right translation by the group has dimension 2, whereas
$\dim(J\cdot W_6)=1$.

A more extreme example, not tabulated, is that of the non-abelian complex
structure $J_0$ on $\Ga\bs W_6$ for which $h^0=3$, $h^1=d=\dim\Kur=6$ and
$\dim\Abel=0$ \cite{Nak,simon}.

\subsection{Final examples}

In general, information on $H^1$ and abelian deformations can be extracted
algebraically using Lemmas~\ref{lemma:closed form} and
Lemma~\ref{def:conditiona}. The computation of $\Kur$ is more challenging,
though it is useful to realize that every parameter is integrable when
$h^1=\dim\Abel$.

In our last two examples, the first applies the theory of \S4.3, whereas the
second replies on this theory to pass directly to a calculation with invariant
differential forms.

\begin{exa} For $(\Ga\bs H_3)\times(\Ga\bs H_3)$ (see Example~\ref{HH}), we
may take \[T_1=\half(X_1-iY_1),\q T_2=\half(X_2-iY_2),\q W=\half(Z_1-iZ_2).\]
The associated complex structural equations are \[[\ol T_1,T_1]=-\half i(W+\ol
W),\qq[\ol T_2,T_2]=\half(W-\ol W).\] In terms of the dual basis $\{\om^i\}$,
any harmonic representative of $H^1$ is a linear combination of \[\oom^1\ot
T_1,\q\oom^2\ot T_2,\q\oom\ot W, \q \oom^1\ot T_2+i\oom^2\ot T_1.\] If
\[\bmu=\mu^1_1\oom^1\ot T_1+\mu^2_2 \oom^2\ot T_2+\mu_3^3\oom\ot W+\mu_1^2
(\oom^1\ot T_2+i\oom^2\ot T_1),\] we obtain \[\bPh=\bmu-\mu_3^3\mu_1^2
(\oom^1\ot T_2-i\oom^2\ot T_1).\] Then (\ref{Abel}) defines an integrable
complex structure.\end{exa}

\def\1{\ol1}\def\2{\ol2}

\begin{exa} The complex structural equations corresponding to Example~\ref{P6}
can be written in the form $d\om^1=0=d\om^2$ and \[2d\om^3=i\om^1\we\oom^1+
\om^1\we\oom^2-\oom^1\we\om^2=i\om^{1\1}+\om^{1\2}-\om^{\12},\] in which the
last expression is an abbreviation of the middle one. By
\cite[Theorem~1.1]{KS}, any invariant complex structure $J'$ sufficiently near
to $J$ has a basis of $(1,0)$ forms that can be written
\begin{equation}\label{basis} \left\{ \begin{array}{l}
\al^1=\om^1+\Phi^1_1\oom^1+\Phi_2^1\oom^2\\
\al^2=\om^2+\Phi^2_1\oom^1+\Phi_2^2\oom^2\\
\al^3=\om^3+\Phi^3_1\oom^1+\Phi_2^3\oom^2+\Phi^3_3\oom^3.\end{array}\right.
\end{equation} This is a dual version of (\ref{Abel}), and the integrabilty
condition (\ref{MC}) amounts to the assertion that $(d\al^3)^{0,2}=0$, or
equivalently \[\begin{array}{rcl} 0 &=& 2d\al^3\we\al^1\we\al^2\\[3pt] &=&
[(i\om^{1\1}+\om^{1\2}-\om^{\12})+\Phi_3^3(i\om^{1\1}-\om^{1\2}+\om^{\12})]
\we[-\Phi^1_2\om^{2\2}+\Phi_2^2\om^{1\2}+\Phi_1^1\om^{\12}]\\[3pt] &=&
\big[-i\Phi^1_2(1+\Phi_3^3)+(1-\Phi_3^3)(\Phi_1^1-\Phi_2^2)\big]\om^{1\12\2}.
\end{array}\] Thus, (\ref{basis}) defines an integrable complex structure on
condition that \[ i(1+\Phi_3^3)\Phi^1_2=(1-\Phi_3^3)(\Phi^1_1-\Phi_2^2),\] and
the coefficients in (\ref{basis}) are sufficiently small (in particular,
$|\Phi_3^3|<1$). 

In this case, the Kuranishi series (\ref{fsum}) is infinite, as it is not
possible to express one coefficient as a polynomial in the others. The term
$\Phi^3_1\oom^1+\Phi^3_2\oom^2$ can be reduced to zero by a suitable right
translation of $J$, and therefore plays no role in the equivalence problem. It
follows that $\dim\Kur= 4$. The abelian condition \[
d\al^3\we\ol\al^1\we\ol\al^2=0\] can be worked out in the same way, and forces
$\Phi^1_2=0$ and $\Phi_1^1=\Phi_2^2$, so $\dim\Abel=3$. \end{exa}

The Table and examples allow us to infer that:

\begin{enumerate}

\item It is possible that every direction in the virtual parameter space is
integrable but only some are tangent to abelian deformation. This occurs for
$(\Ga\bs H_3)\times(\Ga\bs H_3)$ and $\Ga\bs P_6$.

\item It is also possible that some directions are obstructed, irrespective of
the abelian condition. An example is $\Ga\bs W_6$. This phenomenon was
described in \cite[Lemma 4.3]{simon}, and contrasts with the unobstructed
deformation theory for $(\Ga\bs W_6, J_0)$.

\item The centers of $T^6$ and $(\Ga\bs H_3)\times T^3$ certainly have
dimension greater than 1, and these examples do not therefore contradict
Theorem~\ref{fully}.
\end{enumerate}

\n All these observations demonstrate the subtle dependence of $\dim\Kur$ and
$\dim\Abel$ on the underlying algebraic structure of the group $G$.\vs

The techniques of this paper can in theory be applied to study deformations of
the compact quotients of the six series of 2-step nilmanifolds with abelian
complex structures in any complex dimension. 

Other work of the authors shows that in many cases an explicit description of
$\Kur$, and indeed a global moduli space, is possible \cite{GPP,Colin}. It is
also realistic to seek to describe the quotient of the space $\cC(\g)$ by the
group of $\Aut(\g)$ of Lie algebra automorphisms of $\g$, at least near a
generic point of $\cC(\g)$. In this case, in the $W_6$ example,
$\cC(\g)/\hbox{Aut}(\g)$ is locally isomorphic to the quotient of $\Kur$ by
the group of outer automorphisms of $\g$ \cite[\S5]{FPS}.

\parindent0pt\parskip10pt

Department of Mathematics, University of California at Riverside, Riverside, 
CA 92521, USA (\texttt{maclaugh@math.ucr.edu})

Department of Mathematics and Computer Science, University of Southern
Denmark, Campusvej 55, Odense M, DK--5230, Denmark
(\texttt{henrik@imada.sdu.dk})

Department of Mathematics, University of California at Riverside, Riverside, 
CA 92521, USA (\texttt{ypoon@math.ucr.edu})

Mathematics Department, Imperial College, 180 Queen's Gate, London, SW7 2AZ,
UK (\texttt{s.salamon@imperial.ac.uk})

\enddocument